\documentclass{article}
\usepackage[utf8]{inputenc}

\usepackage{graphicx}
\usepackage{tikz}
\usepackage{hyperref}

\title{Enumerating solutions to grid based puzzles with a fixed number of rows}
\author{George Spahn \\
\small{Department of Mathematics, Rutgers University (New Brunswick)}\\
\small{Website: \url{www.georgespahn.com}}}
\date{January 2022}

\begin{document}

\maketitle

In this paper we demonstrate a method for counting the number of solutions to various logic puzzles. Specifically, we remove all of the ``clues" from the puzzle which help the solver to a unique solution, and instead start from an empty grid. We then count the number of ways to fill in this empty grid to a valid solution. We fix the number of rows $r$, vary the number of columns $k$, and then compute the sequence $A_r(k)$, which gives the number of solutions on an empty grid of size $r \times k$.

\section{Ring-Ring}

\includegraphics[scale = 0.5]{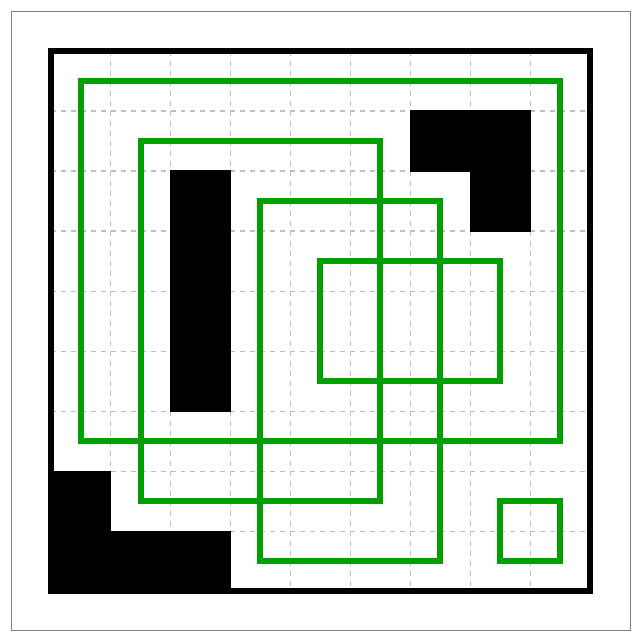}

The New York Times has recently been publishing Ring-Ring puzzles. A solution consists of drawing rectangles so that no grid square remains empty. Rectangles are not allowed to share a side or a corner, however they may overlap. See the diagram for an example of a completed puzzle (the solution is in green). A natural question to ask: How many solutions are there on an empty grid (no black squares) of size $n \times k$? 

\section{}

We consider fixing $k$, and define $A_k(n)$ to be the sequence that gives the number of Ring-Ring solutions on an empty grid of size $n \times k$. This paper gives a method for finding the generating function of this sequence, and gives an explicit formula for the generating function for some small values of $k$.\\
\\
First consider the case $k=4$. Suppose you have a completed solution of an $n \times 4$ grid. Look at one particular column. It is not too difficult to show that it must look like one of the following 15 possibilities. \\
\includegraphics[scale = 0.5]{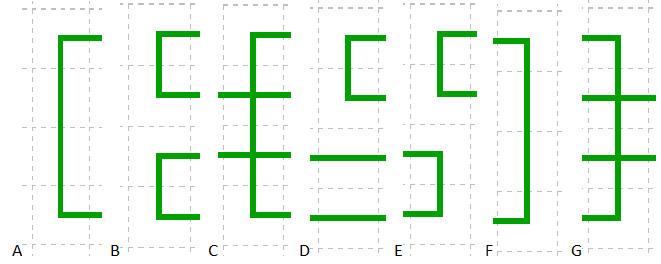}
\includegraphics[scale = 0.5]{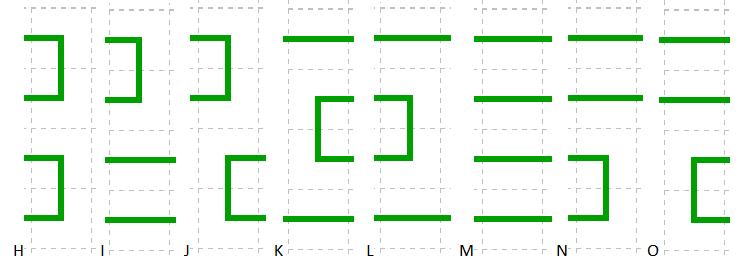}\\
\\
If we call each of the above possibilities a symbol, we have that a solution to the grid must consist of some sequence of symbols. What determines whether such a sequence is legal? We need to ensure that if a symbol leaves some rectangles in progress, that those rectangles are continued in the next symbol. The information of what rectangles are currently in progress can be thought of as a state. When no rectangles are in progress (like is the case for the leftmost column) the only legal symbols are A and B. However if there is currently a rectangle in progress that occupies rows 1 and 2, then the legal symbols are J and O. A state can be described as a set of disjoint subsets of the rows, where each subset is of size 2. These subsets specify the rows in which a rectangle is currently in progress. \\

\begin{center}
\begin{tikzpicture}[scale=0.15]
\tikzstyle{every node}+=[inner sep=0pt]
\draw [black] (37.1,-10.3) circle (3);
\draw [black] (37.1,-10.3) circle (2.4);
\draw [black] (60.1,-10.3) circle (3);
\draw (60.1,-10.3) node {$14$};
\draw [black] (60.1,-27.6) circle (3);
\draw (60.1,-27.6) node {$14,23$};
\draw [black] (60.1,-48.5) circle (3);
\draw (60.1,-48.5) node {$23$};
\draw [black] (18.9,-21.3) circle (3);
\draw (18.9,-21.3) node {$12,34$};
\draw [black] (8.3,-48.5) circle (3);
\draw (8.3,-48.5) node {$12$};
\draw [black] (29.9,-48.5) circle (3);
\draw (29.9,-48.5) node {$34$};
\draw [black] (40.1,-10.3) -- (57.1,-10.3);
\fill [black] (57.1,-10.3) -- (56.3,-9.8) -- (56.3,-10.8);
\draw (48.6,-10.8) node [below] {$A,F$};
\draw [black] (57.1,-10.3) -- (40.1,-10.3);
\fill [black] (40.1,-10.3) -- (40.9,-10.8) -- (40.9,-9.8);
\draw [black] (40.1,-10.3) -- (57.1,-10.3);
\fill [black] (57.1,-10.3) -- (56.3,-9.8) -- (56.3,-10.8);
\draw [black] (40.1,-10.3) -- (57.1,-10.3);
\fill [black] (57.1,-10.3) -- (56.3,-9.8) -- (56.3,-10.8);
\draw [black] (60.1,-13.3) -- (60.1,-24.6);
\fill [black] (60.1,-24.6) -- (60.6,-23.8) -- (59.6,-23.8);
\draw (59.6,-18.95) node [left] {$K,L$};
\draw [black] (57.613,-29.257) arc (-28.59561:-316.59561:2.25);
\draw (52.73,-28.61) node [left] {$M$};
\fill [black] (57.27,-26.63) -- (56.81,-25.81) -- (56.33,-26.69);
\draw [black] (60.1,-24.6) -- (60.1,-13.3);
\fill [black] (60.1,-13.3) -- (59.6,-14.1) -- (60.6,-14.1);
\draw [black] (60.1,-30.6) -- (60.1,-45.5);
\fill [black] (60.1,-45.5) -- (60.6,-44.7) -- (59.6,-44.7);
\draw (59.6,-38.05) node [left] {$C,G$};
\draw [black] (60.1,-45.5) -- (60.1,-30.6);
\fill [black] (60.1,-30.6) -- (59.6,-31.4) -- (60.6,-31.4);
\draw [black] (34.53,-11.85) -- (21.47,-19.75);
\fill [black] (21.47,-19.75) -- (22.41,-19.76) -- (21.89,-18.91);
\draw (25.86,-15.3) node [above] {$B,H$};
\draw [black] (21.47,-19.75) -- (34.53,-11.85);
\fill [black] (34.53,-11.85) -- (33.59,-11.84) -- (34.11,-12.69);
\draw [black] (17.81,-24.1) -- (9.39,-45.7);
\fill [black] (9.39,-45.7) -- (10.15,-45.14) -- (9.21,-44.78);
\draw (12.85,-34.04) node [left] {$N,O$};
\draw [black] (9.39,-45.7) -- (17.81,-24.1);
\fill [black] (17.81,-24.1) -- (17.05,-24.66) -- (17.99,-25.02);
\draw [black] (16.22,-22.623) arc (324:36:2.25);
\draw (11.65,-21.3) node [left] {$M$};
\fill [black] (16.22,-19.98) -- (15.87,-19.1) -- (15.28,-19.91);
\draw [black] (20.02,-24.08) -- (28.78,-45.72);
\fill [black] (28.78,-45.72) -- (28.94,-44.79) -- (28.01,-45.16);
\draw (25.14,-34) node [right] {$D,I$};
\draw [black] (26.9,-48.5) -- (11.3,-48.5);
\fill [black] (11.3,-48.5) -- (12.1,-49) -- (12.1,-48);
\draw (19.1,-48) node [above] {$E,J$};
\draw [black] (11.3,-48.5) -- (26.9,-48.5);
\fill [black] (26.9,-48.5) -- (26.1,-48) -- (26.1,-49);
\draw [black] (28.78,-45.72) -- (20.02,-24.08);
\fill [black] (20.02,-24.08) -- (19.86,-25.01) -- (20.79,-24.64);
\draw [black] (12.8,-7.4) -- (34.12,-9.94);
\draw (12.14,-7.28) node [left] {$START$};
\fill [black] (34.12,-9.94) -- (33.39,-9.35) -- (33.27,-10.35);
\end{tikzpicture}
\end{center}

Above is a finite state machine [S] that computes whether a string of symbols gives a legal solution for $k=4$. Any legal solution must start in the state corresponding to the empty set (no rectangles in progress, labelled with START), and proceed to trace a path in the state machine, eventually returning to the empty state. For example, suppose we proceed from the empty state, to $\{1,4 \}$ and then to $\{1,4 \},\{2,3 \}$ and then to $\{2,3 \}$ and then to $\{1,4 \},\{2,3 \}$ and then to  $\{1,4 \}$ and then to the empty state. This gives the sequence of symbols AKGCLF, which gives the following $6 \times 4$ solution:\\
\includegraphics[scale = 0.7]{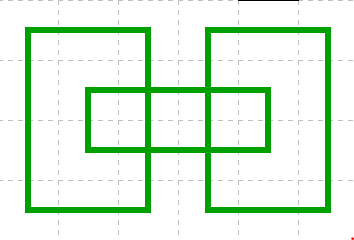}

\section{}
To count the number of solutions on an $n \times 4$ grid, we need to compute the number of paths of length $n$ that both start and end on the empty state. We can solve for the generating function of this sequence directly by solving a linear system of equations. Let $F_i$ be the generating function such that the coefficient of $x^n$ counts the number of paths of length $n$ from the start state that end up at state $i$. Then $$F_i = x \cdot \sum_j F_j  $$ 
where the sum is taken over all states $j$ that have a directed edge to $i$. In words, the number of paths to state $i$ of length $n$ is equal to the sum of the number of paths of length $n-1$ to states that have an edge to state $i$. We now must modify the equation slightly for the start state itself, to account for the fact there is 1 way to get to the start state of length 0. So we have the above equations for $i \neq 0$, and that $$F_0 = 1 + x \cdot \sum_j F_j$$
For $k=4$, we can solve this system of equations to get $$F_0 = \frac{ (1+x)(1-2x)(1-2x-x^2)}{(1-3x-3x^2+10x^3+3x^4-5x^5-x^6) } $$

\section{}
The approach works for any $k$, and the author has maple code which can output the corresponding generating function when given $k$ as input. The states are automatically generated by constructing all possible sets of disjoint subsets of $[k]$, each of size 2, and the edges of the state machine are constructed by ensuring the rules of the puzzle are followed. The corresponding system can be solved automatically by maple. Here is the output for some small values of $k$.\\
\\
$k=5$
$$ \frac{-(2x^2 - 1)(x^4 - 3x^2 + 1)}{(x^8 - 14x^6 + 19x^4 - 8x^2 + 1)} $$
$k=6$
$$\frac{
-(64x^{23} + 518x^{22} - 660x^{21} -\ldots)}{((x + 1)(68x^{24} + 496x^{23} - 1685x^{22} - \ldots))}$$
The growth of these sequences is determined by the largest root of the denominator. We had the computational power to compute this data for $k$ from 2 to 8.

\begin{tabular}{|l|l|}
\hline
2 & 1.618   \\ \hline
3 & 1       \\ \hline
4 & 5.189   \\ \hline
5 & 3.540   \\ \hline
6 & 6.777   \\ \hline
7 & 13.627  \\ \hline
8 & 245.118 \\ \hline
\end{tabular}\\
This trend seems to suggest super exponential growth.

\section{2 Not Touch}
In this section we apply the same technique as the previous section to count the number of 0-1 matrices with fixed column sums and various local constraints. The constraints are inspired by the New York Times puzzles, "2 not touch", in which you must place stars so that they are not adjacent horizontally, vertically, or diagonally, so that the each row and column have exactly 2 stars. Additionally, there must be exactly 2 stars in each region. See the image below for a completed solution. For the purposes of this paper, we remove the 2 stars per region requirement, and count the number of ways to place stars on an empty grid.  We generalize the problem by removing the row sum requirement and letting a(n) count the number of ways to solve it with rows of length n. \\
\\
\includegraphics[scale = 0.7]{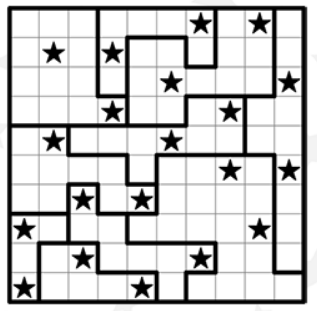}\\
\\
Here the symbols are again the possibilities for each column, and again, we can make a finite state machine to describe which symbols can come after each other in a valid solution. We have code to solve various different versions of the problem, using the method described in the previous section. Let $a_{k,r}(n)$ be the number of $k \times n$ \phantom l 0-1 matrices with the sum of each row equal to $r$. We now add the restrictions:
\begin{itemize}
    \item (i) No two consecutive 1s in a column
    \item (ii) No two consecutive 1s in a row
    \item (iii) No two consecutive 1s diagonally
\end{itemize}
For each subset S of \{(i), (ii),(iii)\} we have code to find the generating function of $a_{k,r,S}(n)$. If S is the empty set, this corresponds to the problem without any of the restrictions.\\
\\
For example, the command gen\_fun(5,2,[1,1,0]) counts the number of solutions with 5 rows, column sum equal to 2, and no two consecutive 1s in a row or column. It produces the output $$ -\frac{2x(x^2 + 6x + 3)}{x^3 + 4x^2 - 1} $$
corresponding to the sequence: 6, 23, 95, 384, 1567, 6361, 25902, 105275, 428363

\section{Circuit Board}
Lastly, we consider the puzzle ``Circuit Board". In this puzzle each dot can be connected with an edge to any subset of its 4 immediate neighbors (up, down, left, right). The goal is to connect all the dots so that each dot has either one edge or three edges emanating from it, without forming any closed loops. See below for an example puzzle. \\
\\
\includegraphics[scale = 0.7]{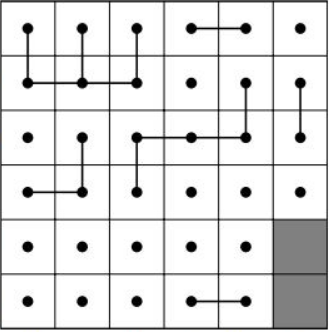}\\
\\
Transitioning to the language of graph theory, we are looking for spanning trees of the grid graph where each vertex has degree 1 or 3. To check whether solutions with $r$ rows are valid, we can easily ensure that each vertex has the appropriate degree as we read the columns. However we must also ensure that the graph is connected and that there are no cycles. To do this, our state also keeps track of which subsets of the current column are connected. Thus a state is essentially a set-partition of the rows, each partition indicating locations of the most recent column that are connected using the columns that have been read so far. \\
\\
The final state consists of set partition with only 1 sub-partition since everything must be connected. In the start state, the set partition is a group of singletons because nothing is connected yet. The number of states is a function of the number of set-partitions on a set with $r$ elements. This number grows super exponentially, and is given by the Bell numbers. The code was able to complete for $r \leq 6$.\\
\\
For $r=4$, it turns out their are never any solutions, regardless of the number of columns. Intrigued by this, I looked for a combinatorial explanation, and discovered that the number of vertices must be congruent to 2 mod 4 for there to be solutions. Any solution can be constructed by starting with 2 vertices that are connected, and then repeatedly expanding degree 1 vertices to be degree 3. If we think of the vertices on a chess board, then some vertices lie on black squares while others lie on white squares. Expanding a vertex either adds 2 black squares or 2 white squares. Thus our final state must have an even number of vertices, and therefore has the same number of white and black squares covered (since it forms a rectangular grid). Thus we must have added the same amount of white squares and black squares, and since we started with 2 squares covered, the total amount of squares covered must be 2 mod 4. \\
\\
For $r=3$ the sequence has the first 10 terms: 0, 1, 0, 0, 0, 10, 0, 0, 0, 36\\
\\ 
The corresponding generating function is $$\frac{x^2(4x^8 - 6x^4 - 1)}{(4x^4 - 1)} $$
Using this, we can show that after the first 10 terms, our sequence satisfies the recurrence $A_3(k) = 4 \cdot A_3(k-4)$. Adding 4 more columns gives a factor of 4 more solutions! \\
\\
For $r=5$ the generating function has approximate degree 50, and for $r=6$ it has approximate degree 100. Download the code to try it yourself.

\section{Code}
My code for the above computations is all available  \href{https://github.com/DarthCalculus/Puzzle-Combinatorics}{\textit{here}}. It is written in maple, and makes use of builtin linear algebra to solve very large systems of equations.

\section{References}
[S] Sipser, M. (2006). Introduction to the theory of computation. Boston: Thomson Course Technology.
\end{document}